# SOLVING *SANGAKU*: A TRADITIONAL SOLUTION TO A NINETEENTH CENTURY JAPANESE TEMPLE PROBLEM

Rosalie Joan Hosking

Abstract. This paper demonstrates how a nineteenth century Japanese votive temple problem known as a *sangaku* from Okayama prefecture can be solved using traditional mathematical methods of the Japanese Edo (1603-1868 CE). We compare a modern solution to a *sangaku* problem from *Sacred Geometry: Japanese Temple Problems* of Tony Rothman and Hidetoshi Fukagawa with a traditional solution of Ōhara Toshiaki (?-1828).

Our investigation into the solution of Ōhara provides an example of traditional Edo period mathematics using the *tenzan jutsu* symbolic manipulation method, as well as producing new insights regarding the contextual nature of the rules of this technique.

## 1. Introduction

During the Japanese Edo period (1603-1868 CE), a unique practice of creating votive tablets with mathematical content emerged. These tablets - known as *sangaku* 算額 - contained a variety of geometrical problems and were dedicated to Shinto shrines and Buddhist temples throughout the country.

This tradition remained largely hidden to historians until the publishing of *Japanese Temple Problems: San Gaku* by David Pedoe and Fukagawa Rothman in 1989 CE. This work provided a collection of problems taken directly from various *sangaku*, along with a selection of solutions using modern mathematical methods and terminology. This was followed up with the more extensive *Sacred Mathematics - Japanese Temple Geometry* by Tony Rothman and Fukagawa Hidetoshi in 2008 CE, which provided further historical and contextual information regarding the *sangaku* tradition.

However, the work of Pedoe, Rothman and Hidetoshi has been criticised for its strictly modern approach to these problems. One critic, Peter J Lu, writes that overall their work "achieves only limited success in showcasing *sangaku* as exemplars of a uniquely Japanese style of mathematics, because that style is never elucidated" [10]. This paper aims to directly respond to such concerns by demonstrating how a nineteenth century Japanese *sangaku* problem solved by Rothman and Hidetoshi with modern methods from the Katayamahiko shrine in Okayama prefecture can be solved using traditional methods of the Japanese Edo period.

The traditional solution used is taken from the work of Ōhara Toshiaki 大原利明 (?-1828), a mathematician from Umeda village in the Adachi ward of modern day Tokyo [6]. In 1810 CE, Ōhara produced the *Sanpō Tenzan Shinan* 算法点竄指南 *Instruction in Tenzan Mathematics*, a three volume work which provided solutions to a number of geometrical problems in mathematics using the traditional manipulation technique known as *tenzan jutsu* 点竄術. The third volume of this work refers to Ōhara as the master of the Ōhara school of mathematics, indicating that he was a mathematical instructor. The method provided by Ōhara is transcribed and translated in section 4.2.

Our investigation into the style of calculation used by Ōhara also provides insight into the contextual nature of calculation in the Edo period, as the use of certain methods appears to have been influenced by the mathematical school one belonged to.

## 2. The Sangaku Tradition

*Sangaku* first appeared in the Edo period, which started when Ieyasu Tokugawa became the military leader in 1603 CE. The new Tokugawa feudal regime instigated many changes, including the national seclusion policy - *sakoku* - which was enforced in 1639 CE. This policy saw scientific works on mathematics, astronomy, and geography become included on lists of banned books [18]. Texts on mathematics reduced in volume due to the ban, prompting the Japanese to begin producing their own works on the subject. This led to the development of a uniquely Japanese mathematical tradition known as *wasan*.

Out of this tradition sprung *sangaku*, where *san* 算 means calculation and *gaku* 額 plaque. These tablets were made from solid wood, and came in a variety of sizes. An examination of 250 *sangaku* listed on the *Sangaku* website of Hiroshi Kotera [9] indicates these tablets had an average width of 153 cm and average height of 69.4 cm. While there are few mentions of the specific type of wood used for *sangaku*, the *Kanchobetsu Kampo Shuroku* refers to a *sangaku* made from *hinoki*, indicating this tree was used for their canvas [7].

*Sangaku* are believed to have been influenced by a pre-existing tradition of hanging tablets in Shino shrines known as *ema* 絵馬 – where *e* 絵 means picture or painting and *ma* 馬 a horse [17]. In ancient times, shrines would ritually sacrifice live horses. However, due to the expense, parishioners began to instead sacrifice paintings of horses on wooden boards as a substitute [8]. During the Muromachi (1333 – 1573 CE) period, larger *ema* began to appear which depicted famous battles, people, or scenes. To house these larger and more artistic tablets, shrines and temples created special areas within their grounds known as *emaden* or *emadō* 絵馬堂 which functioned as informal public art galleries. In the literature, historians such as Shimodaira have described *sangaku* as "mathematical ema" [19] due to their similarities in canvas and location. For example, *sangaku* are commonly found either in the *emaden* of shrines, on the eaves of the shrine or temple halls, or situated attached to a smaller shrine located within the larger complex. *Ema* were also often made from hinoki [23], and could be over one yard in length and width [15].

While *sangaku* may have borrowed their size and location from *ema*, their content was vastly different. *Sangaku* presented mathematical problems which examined the geometry of certain configurations of figures. For example, one might be given a large circle containing five smaller circles of varying sizes and asked to find the value of one of these smaller circles given the diameter of another. These problems ranged in difficulty, and their subject matter included results similar to the Malfatti circles (see [2]), Casey's theorem [3], and Soddy's hexlet theorem [22] which appeared on *sangaku* prior to being known in Europe [17]. These problems were elaborately decorated using pigment paints.

*Sangaku* usually presented their problems in three different sections:

(1) **Problem:** Introductory section describing the diagram and indicating which figure from the diagram the observer needs to find the value of.
(2) **Answer:** Section giving either the numerical value of the sought figure or advising to refer to the next section.
(3) **Formula:** Section which provided a formula for obtaining the solution.

While *sangaku* usually provided a formula for finding the solution, they did not provide any working, leaving the observer to find the solution on their own.

With regard to numbers and labelling, *sangaku* employed the traditional Japanese base ten counting system. The numbers from one to ten in the Japanese style are listed in the right-hand column of Table 1.

| Kanji | Japanese | Reading | Kanji | Japanese | Reading |
|---|---|---|---|---|---|
| 甲 | *kō* | 1st | 一 | *ichi* | one |
| 乙 | *otsu* | 2nd | 二 | *ni* | two |
| 丙 | *hei* | 3rd | 三 | *san* | three |
| 丁 | *tei* | 4th | 四 | *shi/yon* | four |
| 戊 | *bo* | 5th | 五 | *go* | five |
| 己 | *ki* | 6th | 六 | *roku* | six |
| 庚 | *kō* | 7th | 七 | *shichi/nana* | seven |
| 辛 | *shin* | 8th | 八 | *hachi* | eight |
| 壬 | *jin* | 9th | 九 | *kyū* | nine |
| 癸 | *ki* | 10th | 十 | *jyū* | ten |

Table 1. Common Labels Used in Wasan

In this system the largest value is placed before the smaller. For instance, 15 is written *jyu go* 十五, literally 'ten five'. A value such as 20 is expressed *ni jyu* 二十 - 'two ten'. For values in the hundreds and thousands, *hyaku* 百 and *sen* 千 are used in a similar manner.

To label figures and represent unknowns, characters from the Chinese calender such as *kō* 甲, *otsu* 乙, *hei* 丙, and *tei* 丁 (shown in the left-hand column of Table 1) would often be used. These characters were treated as ordinal numbers similar to *a, b, c, d, etc* or first, second, third, etc [4].

## 3. Tenzan Jutsu

One of the major acheivements of Edo period was the development of a symbolic manipulation technique known first as *bōsho hō* 傍書法 and later as *tenzan jutsu*. Sometimes referred to as the method of side-writing [12], it was developed to solve equations which had multiple unknowns.

The *bōsho hō* is accredited to the mathematician Seki Takakazu 関孝和 (? - 1708). The system used vertical lines and kanji characters to represent unknowns and express operations. As the Edo period progressed, the method was later renamed *tenzan jutsu* 点竄術 in the work of the mathematician Matsunaga Yoshinsuke 松永良弼 (1693-1744), who was a second generation pupil of the mathematial school founded by Seki [21].

In *tenzan jutsu*, vertical lines were used to represent numerical values. For example, | 甲 (where as mentioned 甲 *kō* is a character used to represent the first in a series of figures or unknowns. See Table 1) could be used to express 1 甲, || 甲 to express 2 甲, and so on. Negation was shown by a line being crossed through a vertical line. To represent fractions, the numerator would be represented by the kanji character to the right of the line and the denominator by a character on the left, such as 乙 | 甲 for 甲 ÷ 乙. In general, the blank space separating values indicates addition. For example, in Figure 1, the space between | 甲 and 2 | 乙 indicates the calculation | 甲 + 2 | 乙.

In our transcriptions, we translate 甲 as *a*, 乙 as *b*, 丙 as *c*, 丁 as *d*, 戊 as *e*, and so on. Figure 1 displays how *tenzan jutsu* could express the calculation $a + \frac{b}{2} - c$.

$$\begin{array}{ccc} \Big| \, 甲 & \Big| a & (+) \; a \\[1em] 二 \, \Big| \, 乙 & 2 \, \Big| b & (+) \; \dfrac{b}{2} \\[1em] \Big\| \, 小 & \Big\| c & -\; c \end{array}$$

FIGURE 1. Left: Original *tenzan jutsu* calculation. Middle: Transcription. Right: Translation.

**3.1. Rules of Tenzan Jutsu**

In the *Sanpo Tenzan Shinan*, Ōhara includes a set of rules instructing the reader how to manipulate the symbols and lines of the *tenzan jutsu* system to solve a variety of geometrical problems. We have examined and translated the introduction section of the *Sanpo Tenzan Shinan* to produce the following summary of the rules of *tenzan jutsu*.

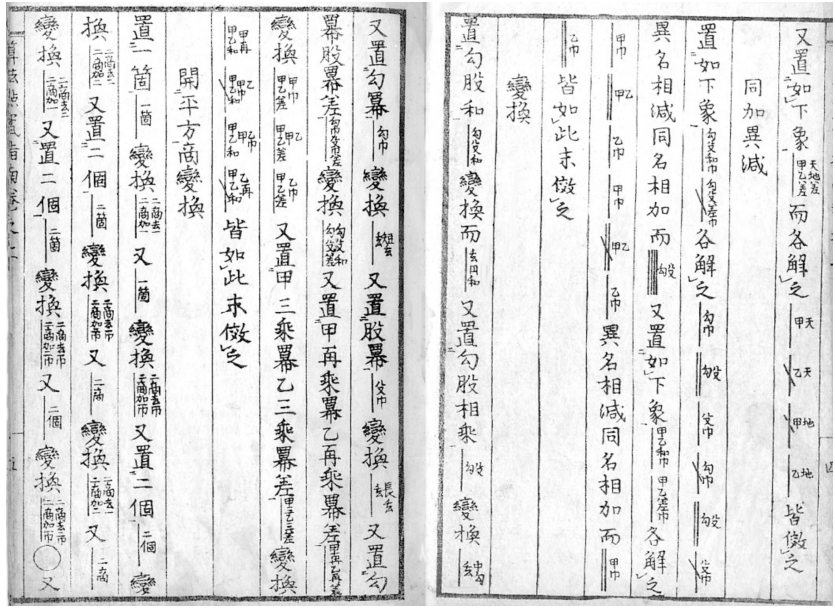

Figure 2. Original text from the *Sanpo Tenzan Shinan* explaining the 'convert' rule

自乘 **Self Multiplication**

The 'self multiplication rule' 自乘 is a means of squaring. For example, applying the rule to $2a$ produces $4a^2$. When there are two terms, squaring and expansion of the square occurs. For example, applying the rule to $a+b$ produces $a^2 + 2ab + b^2$.

括之 **Put Together**

The 'put together' rule 括之 is used to combine elements in various ways. Firstly it can be used to combine terms listed separately with different lines in the *tenzan* system to bring them under one line. For example, say we have the following addition. The rule is also used for combining like terms in the operation equivalent to factorisation.

解之 **Splitting**

The 'splitting' rule has the opposite function of the put together rule. It allows combined terms to be split into individual terms. For instance, if we have $a(b + c)$, the rule changes this into $(a \cdot b) + (a \cdot c)$.

Splitting is also the term used when substitution occurs. When a variable has been defined by a set of terms, the use of the splitting rule indicates that variable is to be decomposed into the terms that make it up. For example, say we have $y + c = x$ and know that $x = a + b$. The splitting rule indicates we can substitute the value of $a + b$ for $x$ and change this expression into $y + c = a + b$.

遍省過乘 **Eliminate Surplus Factors**

When we have an equation such as $3x(a + b + c) - 6ab = 0$, where all terms have a common factor of 3, this rule indicates we are able to eliminate multiplication by three to produce $x(a + b + c) - 2ab = 0$.

同加異減 **Add Same Subtract Different**

This rule states that when we have two terms with similar signs, we can add them, but for different signs we subtract. For example, when there are two similar terms with positive signs - such as $+2ab + 2ab$ - they are combined to form $4ab$. This also applies for two negative values, such that $-2ab - 2ab$ become $-4ab$. However when there are similar terms which have different signs - such as $+a^2 - a^2$ - we subtract. For example, say we have $(a+b)^2 - (a-b)^2$ which produces $a^2 + 2ab + b^2 - a^2 + 2ab - b^2$. These can be grouped into similar terms to form $(+a^2 - a^2)$, $(+2ab + 2ab)$ and $(+b^2 - b^2)$.

變換 **Conversion and** 開平方商變換 **Square Root Converstion**

The 'conversion' rule was used to convert a value into an equivalent value of a different form. In the *Sanpo Tenzan Shinan* the following example, referring to a right-angle triangle, is given:

short side $\cdot$ long side = hypotenuse $\cdot$ altitude

$a^2 - b^2 = (a+b)(a-b)$

The 'square root conversion' rule is used to change a square root into an expanded and equivalent form using binomials. For example, the rule can be used to convert $\sqrt{a}$ into $\sqrt{a}(\sqrt{a} + x)(\sqrt{a} - x)$. As $(\sqrt{a} + x)(\sqrt{a} - x) = 1$, this is equivalent to $\sqrt{a}$. This rule is used when another binomial of either $\sqrt{a} - x$ or $\sqrt{a} + x$ appears in the other terms of the equation. By having the binomial as common to all terms, it can be eliminated from the equation. The following conversions are presented to the reader:

$$1 = (\sqrt{2} - 1)(\sqrt{2} + 1) = (\sqrt{2} - 1)^2(\sqrt{2} + 1)^2 = (\sqrt{3} - 2)(\sqrt{3} + 2)$$
$$= (\sqrt{5} - 2)(\sqrt{5} + 2)$$
$$2 = (\sqrt{2} - 1)^2(\sqrt{2} + 2)^2 = (\sqrt{3} - 1)(\sqrt{3} + 1)$$
$$\sqrt{2} = (\sqrt{2} - 1)(\sqrt{2} + 2) = (\sqrt{2} - 2)^2(\sqrt{2} + 1)^2$$
$$4 = (\sqrt{5} - 1)(\sqrt{5} + 1)$$

The additional conversions are also given, but they are stated to have a positive value when the values are negative:

$$2 = (\sqrt{2} - 2)(\sqrt{2} + 2) = (\sqrt{3} - 2)(\sqrt{3} + 1)^2$$
$$\sqrt{2} = (\sqrt{2} - 2)(\sqrt{2} + 1)$$

乗除括之 **Multiplication and Division Together and** 加減括之 **Addition and Subtraction Together**

With the 'multiplication and division together' rule, when there is one fraction in an equation the rest of the equation is put over a common denominator. It converts an expression such as $\frac{a^2}{b} + 2a + b$ into $\frac{a^2}{b} + \frac{2ab}{b} + \frac{b^2}{b}$.

The rule of 'addition and subtraction together' allows for a negative term to be created in an equation with only positive terms. For instance, say we have $a^2 + ab + b^2$. This then has a negative term added by converting $ab$ into $2ab - ab$ to give $a^2 + 2ab - ab + b^2$. This then becomes $(a+b)^2 - ab$. Another example of this is $a^2 + 2ab + b^2$, which can be converted into $a^2 + 4ab - 2ab + b^2$ and then $(a-b)^2 + 4ab$.

## 4. Katayamahiko Sangaku

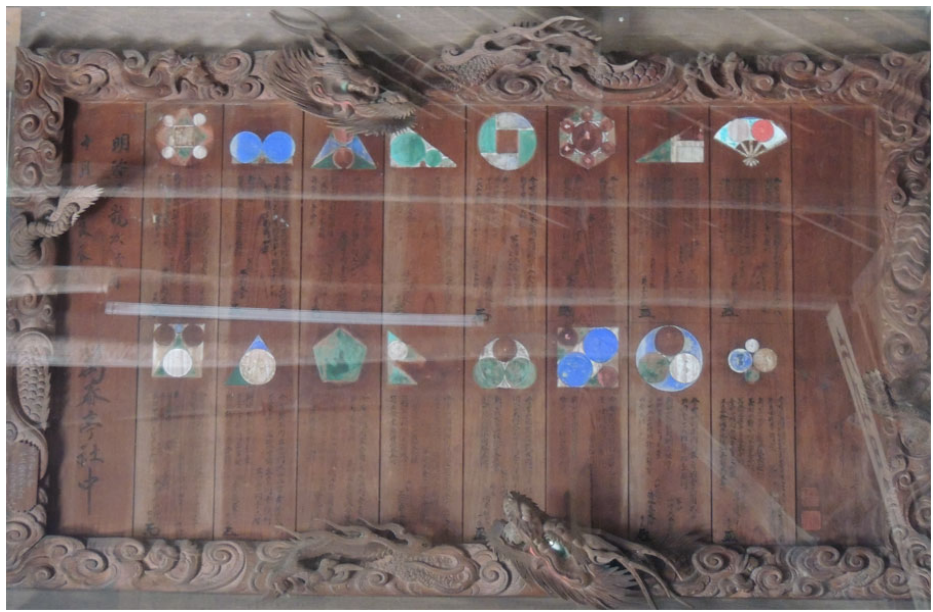

Figure 3. *Sangaku* at Katayamahiko Shrine. (Image by author).

The Katayamahiko *sangaku* was dedicated in October of 1873 CE to the Katayamahiko shrine 片山日子神社 located in the small town of Osafune in Okayama prefecture. It contains sixteen different problems, all of which are beautifully painted and follow the common problem-answer-technique format.

In *Sacred Mathematics - Japanese Temple Geometry*, fifteen of the sixteen problems are described with modern mathematics. The problem we examine is the fourteenth problem, located sixth from the left on the second row on the tablet (see Figure 3). A close up of the diagram is shown in Figure 4. This problem has been chosen because as will be shown, there are interesting features in the traditional solution which differ to modern approaches and provide unique insight into the *wasan* tradition.

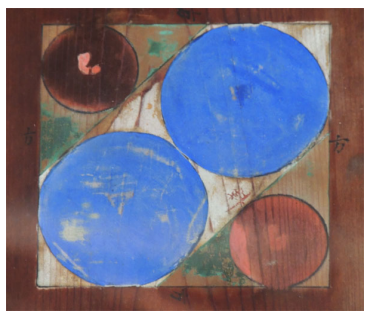

FIGURE 4. The Katayamahiko problem diagram. (Image by author).

The problem text has been translated and transcribed below:

| Transcription | Translation |
| --- | --- |
| 今有如圖方内隔斜容甲乙円各二個只云甲円径一寸問乙円径幾何 | As in the diagram, there is a square segmented by diagonal lines which contains two circles *kō* 甲 and *otsu* 乙. Say the diameter of circle *kō* 甲 is 1 *sun*. Problem - what is the diameter of circle *otsu* 乙? |
| 答云　乙円径五分八厘五毛 | Answer: The diameter of circle *otsu* 乙 is 5 *bu* 8 *rin* 5 *mo* |
| 術日　置二個開平方以域二個余乗甲円径得乙円径合問 | Technique: Put 2 *ko* and take the square root. Subtract 2 *ko*. Multiply the remainder by the diameter of circle *kō* 甲. Obtain the diameter of circle *otsu* 乙 as required. |

### 4.1. Modern Solution

A modern translation of the Katayamahiko problem can be found on page 100 of *Sacred Mathematics - Japanese Temple Geometry*. Rothman and Hidetoshi describe the problem as such:

> [T]wo circles of radius $r$ are inscribed in a square and touch each other at the centre. Each of two smaller circles with radius $t$ touches two sides of the square as well as the common tangent between the two larger circles. Find $t$ in terms of $r$. [17]

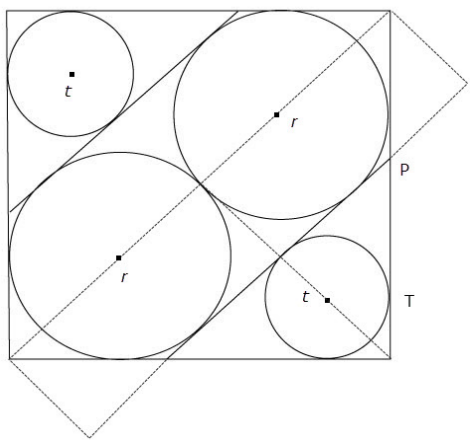

FIGURE 5. Diagram similar to that provided by Rothman and Fukagawa in *Sacred Mathematics - Japanese Temple Geometry* to illustrate their solution to the problem.

They provide a diagram similar to that shown in Figure 5 and the following solution:

> [W]e see that, on the one hand, the length of the central diagonal is $2r+2\sqrt{2}r$. On the other hand, it is also equal to $2PT+2r$. However $PT = t(1+\sqrt{2})$. Equality the two expressions gives:
> $$t = \frac{\sqrt{2}}{\sqrt{2}+1}r = (2-\sqrt{2})r = 0.585786r. \ [17]$$

This solution finds the length of the diagonal of the square in Figure 5 in two different ways. The first is in terms of the larger circles (labelled $r$ in Figure 5) and uses the diagonal which runs from the bottom left-hand corner of Figure 5 to the top right-hand corner. The second is in terms of both the smaller circles (labelled $t$) and larger from the diagonal which runs from the top left-hand corner of the diagram to the bottom right-hand corner. These values - which are equivalent - are then rearranged and used to find the value of the smaller circles in terms of the larger.

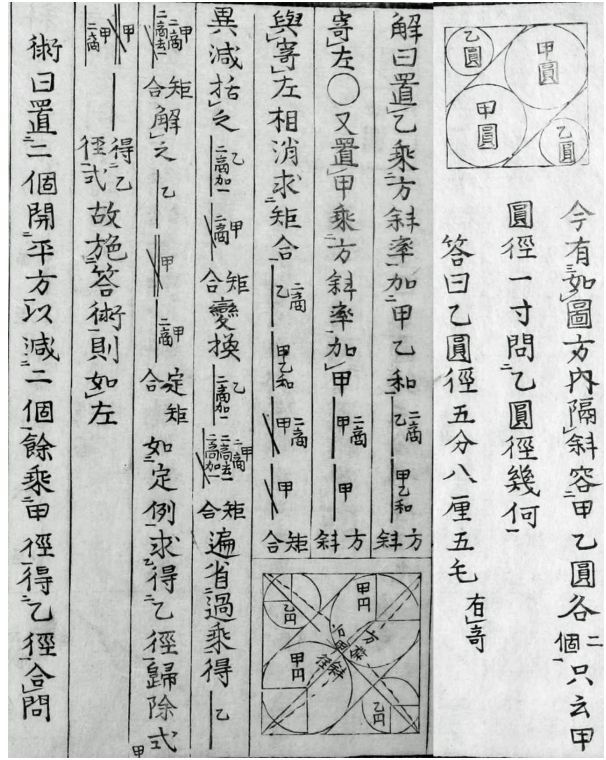

FIGURE 6. Original text from the *Sanpo Tenzan Shinan* showing the calculation using *tenzan jutsu*.

### 4.2. Traditional Solution of Ōhara Toshiaki

As noted, *sangaku*, while providing a formula, did not provide any working, making it difficult for modern historians to ascertain how they were solved. To overcome this issue and locate traditional methods which could be applied to the Katayamahiko problem, we investigated a number of texts created by *wasan* practitioners. This research uncovered the *Sanpo Tenzan Shinan* of Ōhara Toshiaki, which provided instructions on how to solve problems traditionally using the *tenzan jutsu* method. Additionally, this text also provided a solution to a problem almost identical to that found on the Katayamahiko tablet.

In this section, we detail the solution provided by Ōhara. We present a transcription of the original text below in the left column, a translation using modern symbols in the middle column, and translation in modern notation in the right. An image of the original text can be seen in Figure 6. The problem reads:

| Transcription | Translation |
|---|---|
| 今有如圖方内隔斜容甲乙圓各二個只云甲圓径一寸問乙圓径幾何 | As in the diagram, there is a square segmented by diagonal lines which contains two circles *kō* 甲 and *otsu* 乙. Say the diameter of circle *kō* 甲 is 1 *sun*. Problem - what is the diameter of circle *otsu* 乙? |
| 答云　乙圓径五分八厘五毛 | Answer: The diameter of circle *otsu* 乙 is 5 *bu* 8 *rin* 5 *mo* |
| 術曰　置二個開平方以域二個余乘甲圓径得乙圓径合問 | Technique: Put 2 *ko*. Take the square root. Subtract 2 *ko*. Multiply the remainder by the diameter of circle *kō* 甲. Obtain the diameter of circle *otsu* 乙 as required. |

**Given Instructions**:

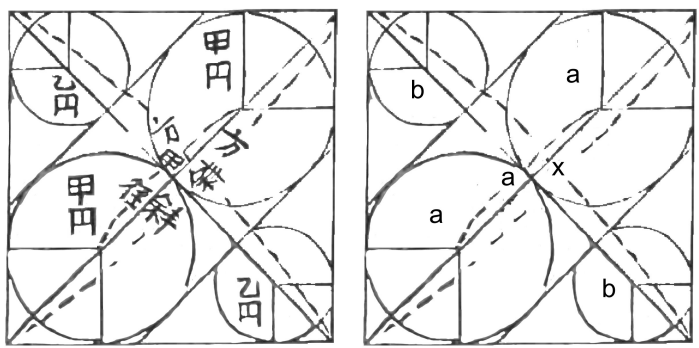

Figure 7. Left: Diagram from the *Sanpo Tenzan Shinan*. Right: Transcription.

[Let 甲 be represented as *a*, 乙 as *b*, and 方斜 (the diagonal) as *x*.]

解曰置乙乘方斜率加甲乙和

Solution: Put 乙 [$b$] and multiply by the diagonal square ratio.[1]. Add 甲 [$a$] plus 乙 [$b$]

| | | |
|---|---|---|
| 乙二<br>商 | $b$ 2<br>√ | $\sqrt{2} \cdot b$ |
| 甲<br>乙<br>和 | $a$<br>$b$<br>$+$ | $(+)\ a + b = x$ |
| 方斜 | x | |

寄左 ◯ 又置甲乘方斜率加甲

Move left. Furthermore, put 甲 [$a$] and multiply by the square diagonal ratio and add 甲 [$a$].

| | | |
|---|---|---|
| 甲二<br>商 | $a$ 2<br>√ | $\sqrt{2} \cdot a$ |
| 甲 | $a$ | $(+)\ a = x$ |
| 方斜 | x | |

與寄左相消求矩合

And move left and cancel to find equation.

| | | |
|---|---|---|
| 乙二<br>商 | $b$ 2<br>√ | $\sqrt{2} \cdot b$ |

---

[1]The diagonal square ratio refers to $\sqrt{2}$.

$$\begin{vmatrix} 甲 \\ 乙 \\ 和 \end{vmatrix} \quad \begin{vmatrix} a \\ b \\ + \end{vmatrix} \quad (+)\ a+b$$

$$\diagdown\!\!\!\begin{vmatrix} 甲 & 二 \\ & 商 \end{vmatrix} \quad \diagdown\!\!\!\begin{vmatrix} a & 2 \\ & \surd \end{vmatrix} \quad -\sqrt{2}\cdot a$$

$$\diagdown\!\!\!\begin{vmatrix} 甲 \end{vmatrix} \quad \diagdown\!\!\!\begin{vmatrix} a \end{vmatrix} \quad -a = 0$$

合矩   0

異域括之
Different terms put together.

$$\begin{vmatrix} 二 & 乙 \\ 商 & \\ 加 & \\ 一 & \end{vmatrix} \quad \begin{vmatrix} 2 & b \\ \surd & \\ + & \\ 1 & \end{vmatrix} \quad b(\sqrt{2}+1)$$

$$\diagdown\!\!\!\begin{vmatrix} 二 & 甲 \\ 商 & \end{vmatrix} \quad \diagdown\!\!\!\begin{vmatrix} 2 & a \\ \surd & \end{vmatrix} \quad -a\cdot\sqrt{2} = 0$$

合矩   0

變換
Convert.

$$\begin{vmatrix} 二 & 乙 \\ 商 & \\ 加 & \\ 一 & \end{vmatrix} \quad \begin{vmatrix} 2 & b \\ \surd & \\ + & \\ 1 & \end{vmatrix} \quad b(\sqrt{2}+1)$$

$$\left\| \begin{matrix} 二二二甲 \\ 商商商 \\ 加去 \\ 一一 \end{matrix} \right. \quad \left\| \begin{matrix} 2 & 2 & 2 & a \\ \sqrt{} & \sqrt{} & \sqrt{} & \\ + & - & & \\ 1 & 1 & & \end{matrix} \right. \quad -a(\sqrt{2}((\sqrt{2}-1)(\sqrt{2}+1)))$$

合矩      0

遍省過乘得
Using elimination of higher power terms many times obtain.

$$\left| \begin{matrix} 乙 \end{matrix} \right| \quad \left| \begin{matrix} b \end{matrix} \right| \quad b$$

$$\left\| \begin{matrix} 二二甲 \\ 商商 \\ 去 \\ 一 \end{matrix} \right. \quad \left\| \begin{matrix} 2 & 2 & a \\ \sqrt{} & \sqrt{} & \\ - & & \\ 1 & & \end{matrix} \right. \quad -a(\sqrt{2}(\sqrt{2}-1)) = 0$$

合矩      0

解之
Splitting.[2]

$$\left| \begin{matrix} 乙 \end{matrix} \right| \quad \left| \begin{matrix} b \end{matrix} \right| \quad b$$

$$\left\| \begin{matrix} 甲 \end{matrix} \right. \quad \left\| \begin{matrix} a \end{matrix} \right. \quad -2a$$

---

[2]We are to separate/expand the terms with a common factor. So $-a(\sqrt{2}(\sqrt{2}-1))$ becomes $-2a + a\sqrt{2}$.

$$\begin{vmatrix} 二甲 \\ 商 \end{vmatrix} \quad \begin{vmatrix} 2 \quad a \\ \surd \end{vmatrix} \quad (+) \; a\sqrt{2} = 0$$

合矩定      0

如定例求得乙径歸除式

As per the typical technique for obtaining the diameter of 乙 [$b$] using the division formula.[3]

$$\begin{vmatrix} 二甲 \\ 商 \end{vmatrix} \; \Big\| \; \begin{vmatrix} 甲 \end{vmatrix} \; \Big\| \; \begin{vmatrix} 2 \; a \\ \surd \end{vmatrix} \; \Big\| \; \begin{vmatrix} a \end{vmatrix} \quad -2a + \sqrt{2}a = b$$

得乙径式      Can find $b$ formula

故絶答術則如左

Therefore absolute answer technique [the answer is] shown left.

## 5. Solution Comparison

It will be observed that both the solution of Fukagawa and Rothman and that of Ōhara rely on finding two different ways of interpreting the diagonal (line $x$ in Figure 6) and then using this information to find the size of the small circles ($b$) in terms of the large ($a$). In *Sacred Mathematics - Japanese Temple Geometry*, where $2r$ from Figure 5 is represented as $a$ and $2t$ as $b$, Rothman and Fukagawa give

(1) $\quad a + \sqrt{2}a = b(\sqrt{2} + 1) + a$

(2) $\quad b = \dfrac{\sqrt{2}}{\sqrt{2} + 1} \cdot a$

However, to get to the end result of $(2 - \sqrt{2})a$ additional steps are required. The denominator must be rationalized, which can be done by multiplying the denominator and numerator by $\sqrt{2} - 1$ as follows

(3) $\quad b = \dfrac{\sqrt{2}(\sqrt{2} - 1)}{(\sqrt{2} + 1)(\sqrt{2} - 1)} \cdot a$

---

[3]What this division formula refers to is not at this point entirely clear.

(4) $\quad b = \dfrac{\sqrt{2}(\sqrt{2}-1)}{1} \cdot a$

(5) $\quad b = (2-\sqrt{2})a$

Though similar, in Ōhara's solution the creation of the fraction in (2) is completely avoided by the use of the 'convert' rule which multiplies $\sqrt{2}a$ by $(\sqrt{2}+1)(\sqrt{2}+1)$ to give

(6) $\quad b(\sqrt{2}+1) - \sqrt{2}a = 0$

(7) $\quad b(\sqrt{2}+1) - \sqrt{2}a(\sqrt{2}+1)(\sqrt{2}-1) = 0$

(8) $\quad b - \sqrt{2}a(\sqrt{2}-1) = 0$

While Ōhara incorrectly rearranges this to give $-2a + \sqrt{2}a = b$, the correct rearranging produces $2a - \sqrt{2}a = b$ which in turn can be simplified as the same $(2-\sqrt{2})a = b$ we find in the solution of Rothman and Fukagawa.

An examination of the problems of the three volumes of the *Sanpo Tenzan Shinan* shows that the 'convert' rule can be used whenever a binomial of the form $\sqrt{a}-x$ or $\sqrt{a}+x$ appears - unless it is in the denominator or part of the final answer. In the *Sanpo Tenshoho Shinan* 算法天生法指南 of Aida Yasuaki 会田安明 (1747–1817), written in 1810 CE, a similar set of rules for *tenzan jutsu* is given to the reader. However, while the lists are very similar, the 'convert' rule of the *Sanpo Tenzan Shinan* is absent from Aida's instructions. This indicates that *tenzan jutsu* did not have a standard set of rules.

As mentioned, Ōhara was the master of the Ōhara school of mathematics, and Aida himself founded a school of mathematics called the Saijyo after falling out with the master of the Seki Takakazu school, Fujita Sadasuke 藤田定資 (1734 - 1807) [20]. As each of these masters presented slightly different rules for calculation with *tenzan jutsu*, it seems that the calculation style used by practitioners was to some degree dependent upon the school they were part of and the texts they had access to.

## 6. Conclusions

While the modern and traditional methods have the same approach to the solution, the calculation steps are different as Ōhara uses the 'convert' rule to avoid creating the fraction $b = \dfrac{\sqrt{2}}{\sqrt{2}+1} \cdot a$ used by Rothman and Hidetoshi.

By examining and comparing these modern and traditional solutions to the Katayamahiko *sangaku* problem, it can be seen that while the modern method provides a correct answer, it does not provide any insight into the original calculation methods of the Japanese and the way in which they would have approached the problem. By investigating an original method, we uncover something interesting about the *wasan* tradition as a whole, as we discover that calculation methods could differ from school to school, meaning *tenzan jutsu* not being standardised.

In conclusion, careful attention must be given to finding and translating original solutions to *sangaku* in the literature in order to accurately represent and understand the methods used in the Edo period.

## Acknowledgement

We would like to thank Tomohiro Uchiyama, Jia-Ming Ying, Clemency Montelle, John Hannah, and Hidetoshi Fukagawa for their help and input with translations and calculations. Sincerest gratitude is also expressed for the Katayamahiko shrine in Osafune for allowing their *sangaku* to be visited and photographed. Lastly the author thanks both the Waseda University Library and Hiroshi Kotera who provided online access to the original mathematical texts used in this research.

## References

[1] Y. Aida, *Sanpo Tenshoho Shinan*, http://www.wasan.jp/archive/tenseihosinan1.pdf, 1810. 会田安明, 算法天生法指南, 1810.
[2] O. Bottema, The Malfatti problem, *Forum Geometricorum* 1 (2001), 43-50.
[3] J. Casey, *A Sequel to Euclid, Hodges*, Figgis, and Co., Dublin, 1881.
[4] N. Dershowitz and E. M. Reingold, *Calendrical Calculations*, Cambridge University Press, New York, 2008.
[5] H. Fukagawa and D. Pedoe, *Japanese Temples Geometry Problems: San Gaku*, Charles Babbage Research Foundation, Winnipeg, 1989.
[6] A. Hirayama, *Gakujutsu o Chūshin to Shita Wasan Shijō no Hitobito*, Tokyo: Fuji Tanki Daigaku Shuppanbu, 1965. 平山諦, 学術を中心とした和算史上の人々, 富士短期大学出版部, 1965.
[7] *Kanchobetsu Kampo Shuroku*, Gihōdō, 1994. 防衛庁防衛年鑑刊行会, 防衛年鑑, 防衛年鑑刊行会, 1994.
[8] E. Kiritani, *Vanishing Japan: Traditions, Crafts and Culture*, Tuttle Publishing, Singapore, 1995.
[9] H. Kotera, Japanese Temple Geometry Problem Sangaku, http://www.wasan.jp.
[10] P. J. Lu, The Blossoming of Japanese Mathematics, *Nature* 454 (2008).
[11] M. Morimoto and T. Ogawa, Mathematical Treatise on Technique of Linkage: An Annotated English Translation of Takebe Katahiro's Tetsujutsu Sankei, *SCIAMVS* 13 (2012), 157-286.
[12] M. Morimoto, The Suanxue Qimen and Its Influence on Japanese Mathematics, *Seki, Founder of Modern Mathematics in Japan. A Commemoration on His Tercentenary*, E. Knobloch, H. Komatsu, and D. Liu (eds), Springer Proceedings in Mathematics and Statistics 39 (2013), 119-132.
[13] T. Ōhara, *Sanpo Tenzan Shinan*, 1810. Waseda University Kotenseki Sogo Database. 大原利明, 算法点竄指南, 1810. 早稲田大古典籍総合データベース.
[14] I. Reader, Letters to the Gods: The Form and Meaning of Ema, *Japanese Journal of Religious Studies* 18(1) (1991), 23-50.
[15] J. Robertson, Ema-gined Community: Votive Tablets (ema) and Strategic Ambivalence in Wartime Japan, *Asian Ethnology* 67(1) (2008), 43-77.
[16] T. Rothman, Japanese Temple Geometry, *Scientific American* 278(5) (1997), 84-91.
[17] T. Rothman and H. Fukagawa, *Sacred Mathematics - Japanese Temple Geometry*, Princeton University Press, Princeton, 2008.
[18] S. Shio, Prohibition of Import of Certain Chinese Books and the Policy of the Edo Government, *Journal of the American Oriental Society* 57(3) (1937), 290-303.

[19] K. Shimodaira, *Mathematics of the Japanese: Wasan*, Kawade Shobō Shinsha, Tokyo, 1972. 下平和夫, 日本人の数学: 和算, 河出書房新社, 1972.
[20] K. Shimodaira, Aida Yasuaki, *Complete Dictionary of Scientific Biography*, ed. Charles Coulsont Gillispie and Frederic Lawrence Holmes and Noretta Koertge and Thomson Gale, Charles Scribner's Sons, Detroit, 2008.
[21] D. E. Smith and Y. Mikami, *A History of Japanese Mathematics*, Cosimo, Inc: New York, 2000.
[22] F. Soddy, The bowl of integers and the hexlet, *Nature* 139 (1927), 77–79.
[23] J. Winter, *East Asian Paintings: Materials, Structures and Deterioration Mechanisms*, Archetype Publications, London, 2008.

Rosalie Hosking
School of Mathematics and Statistics
University of Canterbury
Christchurch, New Zealand
*E-mail address*: `rosalie.hosking@gmail.com`